\documentclass[12pt,twoside]{article}
\usepackage[english]{babel}
\usepackage[latin1]{inputenc}
\usepackage{amsmath}
\usepackage{graphicx}
\usepackage{times,amssymb,amscd}

\numberwithin{equation}{section}

\newcommand{\bA}{\mathbf{A}}

\newcommand{\bE}{\mathbf{E}}

\newcommand{\bH}{\mathbf{H}}

\newcommand{\bL}{\mathbf{L}}

\newcommand{\bR}{\mathbf{R}}
\newcommand{\bS}{\mathbf{S}}
\newcommand{\bV}{\mathbf{V}}

\newcommand{\be}{\mathbf{e}}

\newcommand{\br}{\mathbf{r}}
\newcommand{\bx}{\mathbf{x}}
\newcommand{\by}{\mathbf{y}}

\newcommand{\bT}{\mathbf{T}}

\newcommand{\BV}{\boldsymbol{V}}

\newcommand{\Be}{\boldsymbol{e}}
\newcommand{\Bu}{\boldsymbol{u}}
\newcommand{\Bv}{\boldsymbol{v}}

\newcommand{\cP}{\mathcal{P}}
\newcommand{\cS}{\mathcal{S}}

\newcommand{\EUC}{\mathbf E^3}
\newcommand{\SPH}{\bS^3}
\newcommand{\HYP}{\bH^3}
\newcommand{\SXR}{\bS^2\!\times\!\bR}
\newcommand{\HXR}{\bH^2\!\times\!\bR}
\newcommand{\SLR}{\widetilde{\bS\bL_2\bR}}
\newcommand{\NIL}{\mathbf{Nil}}
\newcommand{\SOL}{\mathbf{Sol}}

\begin{document}
\pagestyle{myheadings}
\markboth{\centerline{Jen\H o Szirmai}}
{On Menelaus' and Ceva's theorems in $\NIL$ geometry}
\title
{On Menelaus' and Ceva's theorems in $\NIL$ geometry
}

\author{Jen\H o Szirmai \\
\normalsize Budapest University of Technology and \\
\normalsize Economics Institute of Mathematics, \\
\normalsize Department of Geometry \\
\normalsize Budapest, P. O. Box: 91, H-1521 \\
\normalsize szirmai@math.bme.hu
\date{\normalsize{\today}}}
\maketitle
\begin{abstract}
In this paper we deal with $\NIL$ geometry, 
which is one of the homogeneous Thurston 3-geometries.
We define the ``surface of a geodesic triangle" using generalized Apollonius surfaces.
Moreover, we show that the ``lines" on the surface of a geodesic triangle can be defined by the famous Menelaus' condition 
and prove that Ceva's theorem for geodesic triangles is true in $\NIL$ space. 
In our work we will use the projective model of $\NIL$ geometry described by E. Moln\'ar in \cite{M97}.

\end{abstract}
\newtheorem{Theorem}{Theorem}[section]
\newtheorem{corollary}[Theorem]{Corollary}
\newtheorem{lemma}[Theorem]{Lemma}
\newtheorem{exmple}[Theorem]{Example}
\newtheorem{definition}[Theorem]{Definition}
\newtheorem{rmrk}[Theorem]{Remark}
\newtheorem{proposition}[Theorem]{Proposition}
\newenvironment{remark}{\begin{rmrk}\normalfont}{\end{rmrk}}
\newenvironment{example}{\begin{exmple}\normalfont}{\end{exmple}}
\newenvironment{acknowledgement}{Acknowledgement}

\section{Introduction}
\label{section1}
In this article, I deal with the generalization and extension of classical 
Euclidean concepts and theorems to Thurston geometries. 
It is very interesting to find theorems that are true in a form in 
all Thurston geometries. 

In our previous paper \cite{Sz21}, we generalized the Menelaus' and Ceva's theorems in $\SXR$ and $\HXR$ spaces
and now we proceed to examine these theorems in $\NIL$ geometry using the notion of the generalized Apollonius surfaces.

The classical definition of the Apollonius circle in the Euclidean plane $\mathbf{E}^2$ is  
the set of all points of $\mathbf{E}^2$ whose distances from two fixed points are in a constant ratio $\lambda\in\mathbf{R}^+$. This definition can be extended
in a natural way to the Thurston geometries 
$$
\EUC,\SPH,\HYP,\SXR,\HXR,\NIL,\SLR,\SOL.
$$
\begin{definition}
The Apollonius surface in the Thurston geometry $X \in \{\EUC,\SPH,\HYP,$ $\SXR,\HXR,\NIL,\SLR,\SOL\}$ is  
the set of all points of $X$ whose geodesic distances from two fixed points are in a constant ratio $\lambda\in\mathbf{R}^+$.
\end{definition}
\begin{rmrk}
A special case of Apollonius surfaces is the geodesic-like bisector (or equidistant) surface ($\lambda=1)$
of two arbitrary points of $X$. These surfaces have an important role in structure of Dirichlet - Voronoi (briefly, D-V) cells and so in the packing and covering problems.
In \cite{PSSz10}, \cite{PSSz11-1}, \cite{PSSz11-2} we studied the geodesic-like 
equidistant surfaces in $\SXR$, $\HXR$ and $\NIL$ geometries, and in \cite{Sz19}, \cite{VSz19}
the translation-like equidistant surfaces in $\SOL$ and $\NIL$ geometries. 
\end{rmrk}

In the present paper, we are interested in {\it geodesic triangles and their surfaces, generalized Menelaus' and Ceva's theorems} in $\NIL$ space \cite{S,T}. 

In Section 2 we describe the projective model and the isometry group of the considered geometry,
moreover, we give an overview about its geodesic curves.

In Section 3 we define the surfaces of geodesic triangles and their {\it lines}, furthermore we examine whether the theorems of Menelaos' and Ceva's are true in $\NIL$ 
space. We prove that the Menelaus' theorem does not follow from the 
structure of $\NIL$ geometry, but plays an important role in defining the surface line, and we show that the Ceva's theorem is true for geodesic triangles in $\NIL$ space.

The computation and the proof is based on the projective model 
of $\NIL$ geometry described by E. Moln\'ar in \cite{M97}.

\section{Projective model of $\NIL$ geometry}
E. {Moln\'ar} has shown in \cite{M97}, that the homogeneous 3-spaces
have a unified interpretation in the projective 3-sphere $\mathcal{PS}^3(\bV^4,\BV_4, \mathbf{R})$. 
In our work we shall use this projective model of $\NIL$ geometry. 
The Cartesian homogeneous coordinate simplex is given by $E_0(\be_0)$,$E_1^{\infty}(\be_1)$,$E_2^{\infty}(\be_2)$,
$E_3^{\infty}(\be_3)$, $(\{\be_i\}\subset \bV^4$ \ and $\text{with the unit point}$ $E(\be = \be_0 + \be_1 + \be_2 + \be_3 ))$. 
Moreover, $\by=c\bx$ with $0<c\in \mathbf{R}$ (or $c\in\mathbf{R}\setminus\{0\})$
defines a point $(\bx)=(\by)$ of the projective 3-sphere $\cP \cS^3$ (or that of the projective space $\cP^3$ where opposite rays
$(\bx)$ and $(-\bx)$ are identified). 
The dual system $\{(\Be^i)\}\subset \BV_4$ describes the simplex planes, especially the plane at infinity 
$(\Be^0)=E_1^{\infty}E_2^{\infty}E_3^{\infty}$, and generally, $\Bv=\Bu\frac{1}{c}$ defines a plane $(\Bu)=(\Bv)$ of $\cP \cS^3$
(or that of $\cP^3$). Thus $0=\bx\Bu=\by\Bv$ defines the incidence of point $(\bx)=(\by)$ and plane
$(\Bu)=(\Bv)$, as $(\bx) \text{I} (\Bu)$ also denote it. Thus, {$\NIL$} can be visualized in the affine 3-space $\bA^3$
(so in $\bE^3$) as well.
\subsection{Geodesic curves and spheres in $\NIL$ space}
In this section we recall the important notions and results from the papers \cite{M97}, \cite{PSSz10}, \cite{Sz13-1}, \cite{Sz07}), \cite{Sz11-1}.

$\NIL$ geometry is a homogeneous 3-space derived from the famous real matrix group $\mathbf{L(R)}$, used by W.~Heisenberg in his electro-magnetic studies.
The Lie theory with the method of projective geometry makes possible to describe this topic.

The left (row-column) multiplication of Heisenberg matrices
     \begin{equation}
     \begin{gathered}
     \begin{pmatrix}
         1&x&z \\
         0&1&y \\
         0&0&1 \\
       \end{pmatrix}
       \begin{pmatrix}
         1&a&c \\
         0&1&b \\
         0&0&1 \\
       \end{pmatrix}
       =\begin{pmatrix}
         1&a+x&c+xb+z \\
         0&1&b+y \\
         0&0&1 \\
       \end{pmatrix}
      \end{gathered} \tag{2.1}
     \end{equation}
defines "translations" $\mathbf{L}(\mathbf{R})= \{(x,y,z): x,~y,~z\in \mathbf{R} \}$ on the points of
$\NIL= \{(a,b,c):a,~b,~c \in \mathbf{R}\}$.
These translations are not commutative, in general. The matrices $\mathbf{K}(z) \vartriangleleft \mathbf{L}$ of the form
     \begin{equation}
     \begin{gathered}
       \mathbf{K}(z) \ni
       \begin{pmatrix}
         1&0&z \\
         0&1&0 \\
         0&0&1 \\
       \end{pmatrix}
       \mapsto (0,0,z)
      \end{gathered}\tag{2.2}
     \end{equation}
constitute the one parametric centre, i.e., each of its elements commutes with all elements of $\mathbf{L}$.
The elements of $\mathbf{K}$ are called {\it fibre translations}. $\NIL$ geometry of the Heisenberg group can be projectively
(affinely) interpreted by the "right translations"
on points as the matrix formula
     \begin{equation}
     \begin{gathered}
       (1;a,b,c) \to (1;a,b,c)
       \begin{pmatrix}
         1&x&y&z \\
         0&1&0&0 \\
         0&0&1&x \\
         0&0&0&1 \\
       \end{pmatrix}
       =(1;x+a,y+b,z+bx+c)
      \end{gathered} \tag{2.3}
     \end{equation}
shows, according to (2.1). Here we consider $\mathbf{L}$ as projective collineation
group with right actions in homogeneous coordinates.

In this context E. Moln\'ar \cite{M97} has derived the well-known infinitesimal arc-length-square, invariant under translations $\bL$ at any point of $\NIL$ as follows
\begin{equation}
   \begin{gathered}
      (dx)^2+(dy)^2+(-xdy+dz)^2=\\
      (dx)^2+(1+x^2)(dy)^2-2x(dy)(dz)+(dz)^2=:(ds)^2
       \end{gathered} \tag{2.4}
     \end{equation}
Hence we get the symmetric metric tensor field $g$ on $\NIL$ by components $g_{ij}$, furthermore, its inverse:
\begin{equation}
   \begin{gathered}
       g_{ij}:=
       \begin{pmatrix}
         1&0&0 \\
         0&1+x^2&-x \\
         0&-x&1 \\
         \end{pmatrix},  \quad  g^{ij}:=
       \begin{pmatrix}
         1&0&0 \\
         0&1&x \\
         0&x&1+x^2 \\
         \end{pmatrix} \\
         \text{with} \ \det(g_{ij})=1.
        \end{gathered} \tag{2.5}
     \end{equation}
The translation group $\mathbf{L}$ defined by formula (2.3) can be extended to a larger group $\mathbf{G}$ of collineations,
preserving the fibering, that will be equivalent to the (orientation preserving) isometry group of $\NIL$.
In \cite{M06} E.~Moln\'ar has shown that
a rotation trough angle $\omega$
about the $z$-axis at the origin, as isometry of $\NIL$, keeping invariant the Riemann
metric everywhere, will be a quadratic mapping in $x,y$ to $z$-image $\overline{z}$ as follows:
     \begin{equation}
     \begin{gathered}
       \br(O,\omega):(1;x,y,z) \to (1;\overline{x},\overline{y},\overline{z}); \\
       \overline{x}=x\cos{\omega}-y\sin{\omega}, \ \ \overline{y}=x\sin{\omega}+y\cos{\omega}, \\
       \overline{z}=z-\frac{1}{2}xy+\frac{1}{4}(x^2-y^2)\sin{2\omega}+\frac{1}{2}xy\cos{2\omega}.
      \end{gathered} \tag{2.6}
     \end{equation}
This rotation formula, however, is conjugate by the quadratic mapping 
     \begin{equation}
     \begin{gathered}
       \mathcal{M}:~x \to x'=x, \ \ y \to y'=y, \ \ z \to z'=z-\frac{1}{2}xy  \ \ \text{to} \\
       (1;x',y',z') \to (1;x',y',z')
       \begin{pmatrix}
         1&0&0&0 \\
         0&\cos{\omega}&\sin{\omega}&0 \\
         0&-\sin{\omega}&\cos{\omega}&0 \\
         0&0&0&1 \\
       \end{pmatrix}
       =(1;x",y",z"), \\
       \text{with} \ \ x" \to \overline{x}=x", \ \ y" \to \overline{y}=y", \ \ z" \to \overline{z}=z"+\frac{1}{2}x"y",
      \end{gathered} \tag{2.7}
     \end{equation}
i.e. to the linear rotation formula. This quadratic conjugacy modifies the $\NIL$ translations in (2.3), as well.
This can also be characterized by the following important classification theorem.
\begin{Theorem}[E.~Moln\'ar \cite{M06} modified]
\begin{enumerate}
\item Any group of $\NIL$ isometries, containing a 3-dimensional translation lattice,
is conjugate by the quadratic mapping in (2.5) to an affine group of the affine (or Euclidean) space $\bA^3=\EUC$
whose projection onto the (x,y) plane is an isometry group of $\bE^2$. Such an affine group preserves a plane
$\to$ point null-polarity.
\item Of course, the involutive line reflection about the $y$ axis
     \begin{equation}
     \begin{gathered}
       (1;x,y,z) \to (1;-x,y,-z),
      \end{gathered} \notag
     \end{equation}
preserving the Riemann metric, and its conjugates by the above isometries in {$1$} (those of the identity component)
are also {$\NIL$}-isometries. There does not exist orientation reversing $\NIL$-isometry.
\end{enumerate}
\end{Theorem}
\begin{rmrk}
We obtain a new projective model of $\NIL$ geometry from the above projective model, derived by the above quadratic mapping $\mathcal{M}$.
This is the {\it linearized model of $\NIL$ space} (see \cite{M06}) that seems to be more advantageous to the future investigations. 
But we remain in the classical so called Heisenberg model in this paper.
\end{rmrk}
\subsection{Geodesic curves and spheres} \label{subsection2}
The geodesic curves of the $\NIL$ geometry are generally defined as having locally minimal arc length between their any two (near enough) points.
The equation systems of the parametrized geodesic curves $g(x(t),y(t),z(t))$ in our model (now by (2.4)) can be determined by the
Levy-Civita theory of Riemann geometry.
We can assume, that the starting point of a geodesic curve is the origin because we can transform a curve into an
arbitrary starting point by translation (2.1);
\begin{equation}
\begin{gathered}
        x(0)=y(0)=z(0)=0; \ \ \dot{x}(0)=c \cos{\alpha}, \ \dot{y}(0)=c \sin{\alpha}, \\ \dot{z}(0)=w; \ - \pi \leq \alpha \leq \pi. \notag
\end{gathered}
\end{equation}
The arc length parameter $s$ is introduced by
\begin{equation}
 s=\sqrt{c^2+w^2} \cdot t, \ \text{where} \ w=\sin{\theta}, \ c=\cos{\theta}, \ -\frac{\pi}{2}\le \theta \le \frac{\pi}{2}, \notag
\end{equation}
i.e. unit velocity can be assumed.

The equation systems of a helix-like geodesic curves $g(x(t),y(t),z(t))$ if $0<|w| <1 $:
\begin{equation}
\begin{gathered}
x(t)=\frac{2c}{w} \sin{\frac{wt}{2}}\cos\Big( \frac{wt}{2}+\alpha \Big),\ \
y(t)=\frac{2c}{w} \sin{\frac{wt}{2}}\sin\Big( \frac{wt}{2}+\alpha \Big), \notag \\
z(t)=wt\cdot \Big\{1+\frac{c^2}{2w^2} \Big[ \Big(1-\frac{\sin(2wt+2\alpha)-\sin{2\alpha}}{2wt}\Big)+ \\
+\Big(1-\frac{\sin(2wt)}{wt}\Big)-\Big(1-\frac{\sin(wt+2\alpha)-\sin{2\alpha}}{2wt}\Big)\Big]\Big\} = \\
=wt\cdot \Big\{1+\frac{c^2}{2w^2} \Big[ \Big(1-\frac{\sin(wt)}{wt}\Big)
+\Big(\frac{1-\cos(2wt)}{wt}\Big) \sin(wt+2\alpha)\Big]\Big\}. \tag{2.8}
\end{gathered}
\end{equation}
In the cases $w=0$ the geodesic curve is the following:
\begin{equation}
x(t)=c\cdot t \cos{\alpha}, \ \ y(t)=c\cdot t \sin{\alpha}, \ \ z(t)=\frac{1}{2} ~ c^2 \cdot t^2 \cos{\alpha} \sin{\alpha}. \tag{2.9}
\end{equation}
The cases $|w|=1$ are trivial: $(x,y)=(0,0), \ z=w \cdot t$.
\begin{definition}
The distance $d(P_1,P_2)$ between the points $P_1$ and $P_2$ is defined by the arc length of geodesic curve
from $P_1$ to $P_2$.
\end{definition}
\begin{definition}
 The geodesic sphere of radius $R$ with centre at the point $P_1$ is defined as the set of all points 
 $P_2$ in the space with the condition $d(P_1,P_2)=R$. Moreover, we require that the geodesic sphere is a simply connected 
 surface without selfintersection 
 in the $\NIL$ space.
 \begin{rmrk}
 We shall see that this last condition depends on radius $R$.
 \end{rmrk}
 \end{definition}
 \begin{definition}
 The body of the geodesic sphere of centre $P_1$ and of radius $R$ in the $\NIL$ space is called geodesic ball, denoted by $B_{P_1}(R)$,
 i.e., $Q \in B_{P_1}(R)$ iff $0 \leq d(P_1,Q) \leq R$.
 \end{definition}
 We proved in\cite{Sz13-2,Sz07} the following important theorems:
 \begin{Theorem}[\cite{Sz07}]
 The geodesic sphere and ball of radius $R$ exists in the $\NIL$ space if and only if $R \in [0,2\pi].$
 \end{Theorem}
 \begin{Theorem}[\cite{Sz13-2}]
 The geodesic $\NIL$ ball $B(S(R))$ is convex in affine-Euclidean sense in our model if and only if $R \in [0,\frac{\pi}{2}]$. 
\end{Theorem}
\subsection{Some properties of geodesic curves and spheres }
In the following, we determine some important properties of geodesic curves and spheres, which we will use in the following sections. 

\begin{enumerate}
\item Consider points $P(x(t),y(t),z(t))$ lying on a sphere $S$ of radius $R$  
centred at the origin. The coordinates of $P$ are given by parameters $(\alpha\in[-\pi, \pi), ~ \theta\in [-\frac{\pi}{2},\frac{\pi}{2}], ~ R>0)$ (see (2.8), (2.9)). 

We obtain directly from the equations (2.8) and (2.9) the following

\begin{lemma}
\begin{enumerate}
\item
$$x(t)^2+y(t)^2=\frac{4c^2}{w^2}\sin^2{\frac{wt}{2}},$$
that means, that if $\theta \ne \pm \frac{\pi}{2}$ and $t=R$ is given and $\alpha\in[-\pi,\pi)$ then the endpoints $P$ of the geodesic curves lie on a cylinder of radius 
$r=\left| \frac{4c}{w}\sin{\frac{wR}{2}}\right|$ with axis $z$. Therefore, we obtain the following connection between parameters $\theta$ and $R$:
\begin{equation}
R=2\cdot \arcsin \left[  \frac{\sqrt{x^2(R)+y^2(R)}}{2\cdot \cot{\theta}}\right] \frac{1}{\sin{\theta}} \tag{2.10}
\end{equation}
\item
If $\theta = \pm \frac{\pi}{2}$ then the endpoints $P(x(R),y(R),z(R))$ of the geodesics $g(x(t),y(t),z(t))$ lie on the 
$z$-axis thus their orthogonal projections onto the $[x,y]$-plane is the origin and $x(R)=y(R)=0$, $z(R)=d(O,R)=R$.

\item
Moreover, the cross section 
of the spheres $S$ with the plane $[x,z]$ is given by the following system of equation:
     \begin{equation}
     \begin{gathered}
    X(R,\theta)=\frac{2c}{w} \sin{\frac{wR}{2}}=\frac{2\cos{\theta}}{\sin{\theta}} \sin{\frac{R \sin{\theta}}{2}}, \\ 
    Z(R,\theta)=wR+\frac{c^2R}{2w} - \frac{c^2}{2w^2}\sin{wR}= \\ 
    R\sin{\theta}+\frac{R\cos^2{\theta}}{2\sin{\theta}} - \frac{\cos^2{\theta}}{2\sin^2{\theta}}\sin(R\sin{\theta}), \ \ (\theta\in [-\frac{\pi}{2},\frac{\pi}{2}]\setminus\{0\}); \\
    \text{if} \ \theta=0 \ \text{then} \ X(R,0)=R, \ Z(R,0)=0. \tag{2.11}
   \end{gathered}
   \end{equation}
\end{enumerate}
\end{lemma}
\begin{rmrk}
The parametric equations of the geodesic sphere of radius $R$ can be generated from (2.1) by $\NIL$ rotation (see (2.6)).
\end{rmrk}
\item 
We introduce the usual notion of the fibre projection $\mathcal{P}$ that is a projection parallel to fibre lines (parallel to $z$-axis), onto the $[x,y]$ plane.
The image of a point $P$ is the intersection with the $[x,y]$ base plane of the line parallel to fibre line passing through $P$, $\mathcal{P}(P)=P^*$.

Analysed the parametric equations of the geodesic curves $g(x(t),y(t),z(t))$ with starting points at the origin we get the following
\begin{lemma}
If $0 <|w| <1$ for geodesic curve $g(x(t),y(t),z(t))$ $(t\in [0,R])$ then 
the fibre projection $\mathcal{P}$ of the geodesic curves onto the $[x, y]$ plane is an Euclidean circle arc where it is contained by a circle with equation 
\begin{equation}
\Big(x(t)+\frac{c}{w}\sin{\alpha}\Big)^2+\Big(y(t)-\frac{c}{w}\cos{\alpha}\Big)^2=\Big(\frac{c}{w}\Big)^2=\cot^2{\theta}. \tag{2.12}
\end{equation}
If $w=0$ then fibre projection $\mathcal{P}$ of the geodesic curves $g(x(t),y(t),z(t))$ $(t\in [0,R])$ onto the $[x, y]$ plane is a segment with starting point at the origin
where it is contained by the straight line with equation
\begin{equation}
y=\tan{\alpha}\cdot x.\tag{2.13}
\end{equation}
If $w=1$ then the fibre projection $\mathcal{P}$ of the geodetic curves $g(x(t),y(t),z(t))$ $(t\in [0,R])$ onto 
the $[x, y]$ plane is the origin.
\end{lemma}
We get directly from the equation (2.12) the following
\begin{corollary}
\begin{enumerate}
\item If we know the equation of the circle that contains the orthogonal projected image $OP^*$ of a geodesic curve segment 
$g_{OP}=g(x(t),$ $y(t),z(t))$ $(t\in [0,R])$ onto the $[x,y]$ plane where $0 <|w| <1$ is a known real number and the 
coordinates of $P^*=(x(R), y(R),0)$ then the parametric equation of the geodesic curve segment $g_{OP}$ is uniquely determined. 
That means that there is {\it one-to-one correspondence between the circle arcs $OP^*$ and the geodesic curve segments $OP$} by the above sense. 
\item If $w=0$ then the fibre projection $\mathcal{P}$ of the geodetic curves is a segment with starting point at the origin
where it is contained by the straight line $y=\tan{\alpha}\cdot x$, therefore in this situation it is 
{\it one-to-one correspondence between the projected image $OP^*$ and the geodesic curve segments $OP$}, too. 
\item If $w=1$ then the fibre projection $\mathcal{P}$ of the geodetic curves 
is the origin so here it is also a {\it one-to-one correspondence between the projected image and the above geodesic curves.}
\end{enumerate}
\end{corollary}
\end{enumerate}
\section{Geodesic triangles and their surfaces}
We consider $3$ points $A_0$, $A_1$, $A_2$ in the Heisenberg model of $\NIL$ space (see Section 2).
The {\it geodesic segments} $a_k$ between the points $A_i$ and $A_j$
$(i<j,~i,j,k \in \{0,1,2\}, k \ne i,j$) are called sides of the {\it geodesic triangle} with vertices $A_0$, $A_1$, $A_2$.
It can be assumed by the homogeneity of the $\NIL$ geometry that $A_0=(1,0,0,0)$.
\begin{figure}[ht]
\centering
\includegraphics[width=12cm]{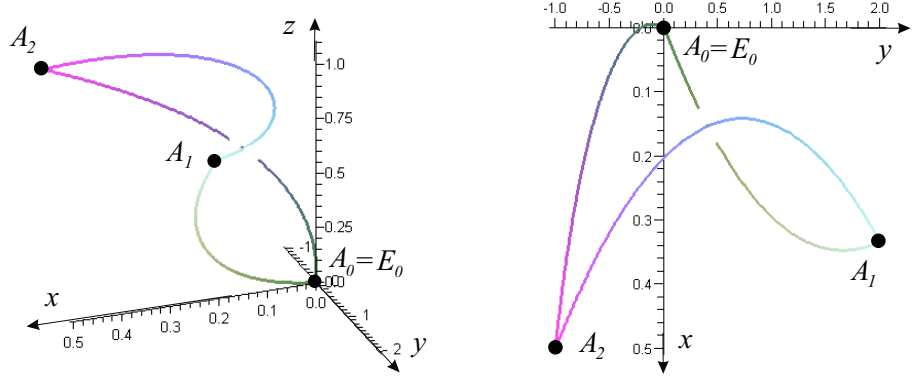}
\caption{Two different views of geodesic triangle with vertices $A_1=(1,0,0,0)$, $A_2=(1,1/2,-1,1)$, $A_3=(1,1/3,2,1)$.}
\label{}
\end{figure}
However, defining the surface of a geodesic triangle in $\NIL$ space is not straightforward. The usual geodesic triangle surface definition is 
not possible because the geodesic curves starting from different vertices and ending at points of the corresponding opposite edges define different 
surfaces, i.e. {\it geodesics starting from different vertices and ending at points on the corresponding opposite side usually do not intersect 
(see Fig.~2).}
\begin{figure}[ht]
\centering
\includegraphics[width=12cm]{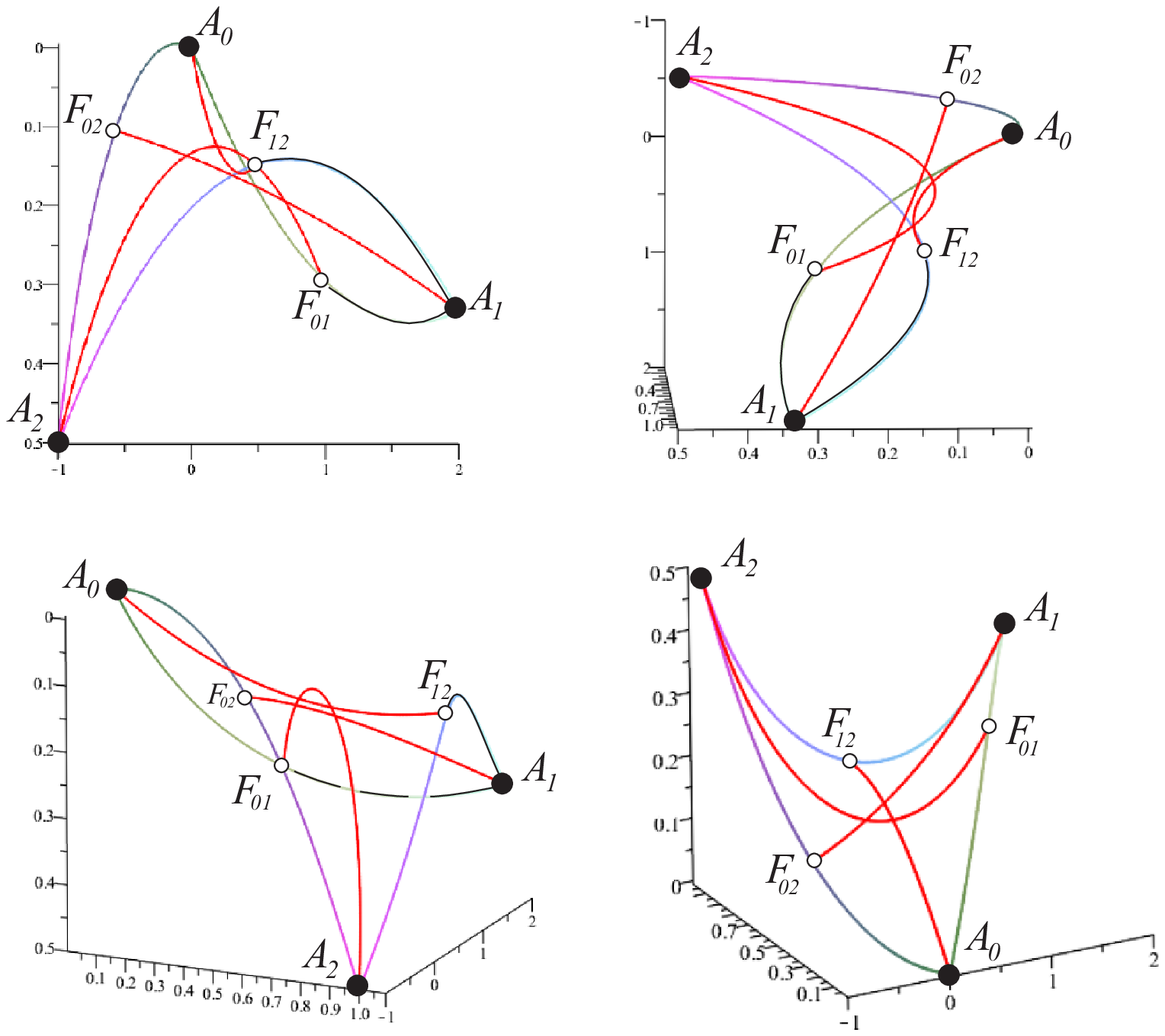}
\caption{$A_0A_1A_2$ is a geodesic triangle in $\NIL$ space with vertices $A_0=(1,1,0,0)$, $A_1=(1,1/3,2,1)$, $A_2=(1,1/2,-1,1)$, $F_{ij}$ is the midpoint of 
the geodesic segment $A_iA_j$ $(i,j\in \{0,1,2\}, j>i)$. The figures show that the geodesic segments $A_kF_{ij}$ $(k\in \{0,1,2\}, k\ne i,j)$ do not intersect each other. 
}
\label{}
\end{figure}
We use for the definition of surfaces of geodesic triangles similarly to the $\SXR$ and $\HXR$ spaces (see \cite{Sz21}) the generalization of the Apollonius surfaces.

The extension of the classical definition of the Apollonius circle of the Euclidean plane $\mathbf{E}^2$ to Thurston geometries is the following  
\begin{definition}
The Apollonius surface $\mathcal{A}\cS^X_{P_1P_2}(\lambda)$ in the Thurston geometry $X$ is  
the set of all points of $X$ whose geodesic distances from two fixed points are in a constant ratio $\lambda\in\mathbf{R}^+_0$ where $X\in \{\EUC,\SPH,\HYP,\SXR,\HXR,\NIL,\SLR,\SOL\}.
$ i.e. $\mathcal{A}\cS^X_{P_1P_2}(\lambda)$ of two arbitrary points $P_1,P_2 \in X$ consists of all points $P'\in X$,
for which $d^X(P_1,P')=\lambda \cdot d^X(P',P_2)$ ($\lambda\in [0,\infty$) where $d^X$ is the corresponding distance function of $X$. If $\lambda=0$, 
then $\mathcal{A}\cS^X_{P_1P_2}(0):=P_1$ and it is clear, that in case $\lambda \to \infty$ then $d(P',P_2) \to 0$ therefore we say $\mathcal{A}\cS^X_{P_1P_2}(\infty):=P_2$.
\end{definition}
We introduce a new definition of the surface $\mathcal{S}_{A_0A_1A_2}$ of the geodesic triangle by the following steps:
\begin{definition}
\begin{enumerate}
\item We consider the geodesic triangle $A_0A_1A_2$ in the projective model of $\NIL$ space and consider the Apollonius surfaces 
$\mathcal{A}\cS_{A_0A_1}$ $(\lambda_1)$ and $\mathcal{A}\cS_{A_2A_0}(\lambda_2)$ ($\lambda_1,\lambda_2 \in [0,\infty)$, $\lambda_1^2+\lambda_2^2>0$). It is clear, that if 
$Y \in \mathcal{C}(\lambda_1,\lambda_2):=\mathcal{A}\cS_{A_0A_1}(\lambda_1)\cap \mathcal{A}\cS_{A_2A_0}(\lambda_2)$ then 
$\frac{d(A_0,Y)}{d(Y,A_1)}=\lambda_1$ and  $\frac{d(A_2,Y)}{d(Y,A_0)}=\lambda_2$ $\Rightarrow$ $\frac{d(A_2,Y)}{d^X(Y,A_1)}=\lambda_1 \cdot \lambda_2$
for parameters $\lambda_1,\lambda_2 \in (0,\infty)$ and if $\lambda_1=0$ then $\mathcal{C}(\lambda_1,\lambda_2)=A_0$, if $\lambda_2=0$ then $\mathcal{C}(\lambda_1,\lambda_2)=A_2$ 
\item 
\begin{equation}
\begin{gathered}
P(\lambda_1,\lambda_2):=\{ P \in \NIL~|~P \in \mathcal{C}(\lambda_1,\lambda_2)~{\text{and}}~ d(P,A_0)=\min_{Q \in \mathcal{C}(\lambda_1,\lambda_2)}({d(Q,A_0)}) \\
~\text{with given real parameters}~\lambda_1,\lambda_2 \in [0,\infty),~ \lambda_1^2 + \lambda_2^2 > 0 \}
\end{gathered} \tag{3.10}
\end{equation}
\item The surface $\mathcal{S}_{A_0A_1A_2}$ of the geodesic triangle $A_0A_1A_2$ is 
\begin{equation}
\mathcal{S}_{A_0A_1A_2}:=\{P(\lambda_1,\lambda_2) \in \NIL,~ \text{where}~\lambda_1,\lambda_2 \in [0,\infty),~ \lambda_1^2 + \lambda_2^2 > 0 \}. \tag{3.11}
\end{equation}
\end{enumerate}
\end{definition}

We introduce the following notations: 

{\it 
1. If the vertices of the geodesic triangle are contained by a Euclidean plane parallel to fibre lines  (parallel to the $z$ axis) 
then the triangle is called fibre type triangle. 

2. In other cases the triangle is in general type.}
\section{On Menelaus' and Ceva's theorems in $\NIL$ space} 
First we recall the usual definition of simply ratios in Euclidean plane $\bE^2$ plane. 
\begin{definition}
If $A$, $B$ and $P$ are distinct points on a line in the Euclidean plane $\bE^2$, then 
their simply ratio is
$s^E(A,P,B) =  d^E(A,P))/d^E(P,B))$, if $P$ is between $A$ and $B$, and
$s^E(A,P,B) = -d^E(A,P))/d^E(P,B))$, otherwise where
$d^E$ denotes the Euclidean distance function. 
\end{definition}

Note that the value of $s^E(A,P,B)$ determines the position of $P$ relative to $A$ and $B$.

\begin{Theorem} [Menelaus' theorem for triangles in $\bE^2$ plane]
If is a $l$ line not through any vertex of a triangle $ABC$ such that
$l$ meets $BC$ in $Q$, $AC$ in $R$, and $AB$ in $P$,
then $$s^E(A,P,B)s^E(B,Q,C)s^E(C,R,A) = -1.~ ~ \square $$
\end{Theorem}

\begin{Theorem}[Ceva's theorem for triangles in $\bE^2$ plane]
If $T$ is a point not on any side of a triangle $ABC$ such that
$AT$ and $BC$ meet in $Q$, $BT$ and $AC$ in $R$, and $CT$ and $AB$ in $P$,
then $$s^E(A,P,B)s^E(B,Q,C)s^E(C,R,A) = 1. ~ ~ \square $$
\end{Theorem}
\begin{rmrk}
The converses of Menelaus' and Ceva's theorems are also true.  
\end{rmrk}
\subsection{How to define the notion of {\it lines} on the surfaces of geodesic triangles?}
If we want to discuss on Menelaus' and Ceva's theorems, 
we must first clarify what we consider to be a {\it line} on the surface of a geodesic triangle and what is the definition of the simple ratio.

Let $\mathcal{S}_{A_0A_1A_2}$ be the surface of the geodesic triangle $A_0A_1A_2$ and $P_1$, $P_2 \in \mathcal{S}_{A_0A_1A_2}$ two given points. 
Natural requirements for a {\it line} through points $P_1$ and $P_2$ lying on $\mathcal{S}_{A_0A_1A_2}$:
\begin{enumerate}
\item Two surface points uniquely determine one {\it line} (connecting curve) $\mathcal{G}^{\mathcal{S}_{A_0A_1A_2}}_{P_1P_2}$.
\item Any two points on a surface {\it line} $\mathcal{G}^{\mathcal{S}_{A_0A_1A_2}}_{P_1P_2}$ define the same one.
\item The surface {\it line} determined by two points of a geodesic curve lying on the surface $\mathcal{S}_{A_0A_1A_2}$ coincide with the geodesic curve.   
\end{enumerate}
\begin{rmrk}
An obvious option for definition of a line (connecting curve) $\mathcal{G}^{\mathcal{S}_{A_0A_1A_2}}_{P_1P_2}$ would be the fibre projection of the geodesic curve $g_{P_1P_2}$ 
into the surface $\mathcal{S}_{A_0A_1A_2}$ but it is clear, that this definition does not satisfy the requirement 2. 
\end{rmrk}
We consider a {\it geodesic triangle $A_0A_1A_2$} in the projective model of $\NIL$ space (see Section 3).
Without limiting generality, we can assume that $A_0=(1,0,0,0)$. 
The geodesic lines that contain the sides $A_0A_1$ and $A_0A_2$ of the given triangle can be characterized directly 
by the corresponding parameters $\theta_i$ and $\alpha_i$ $(i=1,2)$
(see (2.8) and (2.9)).
The geodesic curve including the side segment $A_1A_2$ is also determined by one of its endpoints and its parameters. 
In order to determine the corresponding parameters of this 
geodesic line we use {\it e.g. a $\NIL$ translation} $\bT(A_1)$, as elements of the isometry group of $\NIL$ geometry, that
maps the $A_1=(1,x_1,y_1,z_1)$ onto $A_0=(1,0,0,0)$ (up to a positive determinant factor). 

\begin{rmrk} 
Because of the results of Theorem 2.4, we assume, that the surface $\mathcal{S}_{A_0A_1A_2}$ of the geodesic triangle 
$A_0A_1A_2$ is contained in a geodesic 
$\NIL$ sphere of radius $\pi$. 
\end{rmrk}

First we generalize the notion of simple ratio to the point triples lying on geodesic lines of the $\NIL$ space:
\begin{definition}
If $A$, $B$ and $P$ be distinct points on a geodesic curve in the $\NIL$ space then 
their simply ratio is
$$s^N(A,P,B) =  {d(A,P)}/{d(P,B)},$$ if $P$ is between $A$ and $B$, and
$$s^N(A,P,B) = -{d(A,P)}/{d(P,B)},$$ 
otherwise, where
$d$ is the distance function of $\NIL$ geometry.   
\end{definition}
Let $A$, $B$ and $P$ be distinct points on a non-fibrum-like geodetic curve in the $\NIL$ and let $A^*$, $B^*$ and $P^*$ 
their projected images by $\mathcal{P}$. 
The geodesic curve segment $g_{AB}$ also determined by parameters 
$(t,\theta,\alpha)$ (see (2.8), (2.9)) 
and the parameters of geodesic curve segment $g_{AP}$
are $(t_p,\theta_p=\theta,\alpha_p=\alpha)$. Their images $\overset{\LARGE\frown}{A^*P^*}$ and $\overset{\LARGE\frown}{A^*B^*}$ by fibre projection are circle arcs or/and line segments that are determined 
by parameters $(t_p\cdot \cos{\theta}, \theta,\alpha)$ and $(t \cdot \cos{\theta}, \theta, \alpha)$ (see Lemma 2.6, Lemma 2.8 and Corollary 2.9).

\begin{lemma}
The Euclidean length $\mathcal{C}(A^*,P^*)$ of circle arc or line segment $\overset{\LARGE\frown}{A^*P^*}$ satisfy the following equations 
\begin{equation}
\begin{gathered}
\mathcal{C}(A^*,P^*)=d(A,P)\cdot \cos{\theta},\ \ \ \mathcal{C}(P^*,B^*)=d(P,B)\cdot \cos{\theta}.
\end{gathered} \tag{4.1}
\end{equation}
{Therefore, the projection $\mathcal{P}$  preserves the ratio of lengths by the above sense.}
\end{lemma}

\begin{definition}
Let $\mathcal{S}_{A_0A_1A_2}$ be the surface of the geodesic triangle $A_0A_1A_2$ and $P_1$, $P_2 \in \mathcal{S}_{A_0A_1A_2}$ two given point. 

1. If the points $P_1$ and $P_2$ lie on a fibre line then the connecting curve 
$\mathcal{G}^{\mathcal{S}_{A_0A_1A_2}}_{P_1P_2} \subset \mathcal{S}_{A_0A_1A_2}$ is coincides with the $P_1P_2$ segment on the fibre line. 

2. In other cases the {\it line} (connecting curve) $\mathcal{G}^{\mathcal{S}_{A_0A_1A_2}}_{P_1P_2}$ is the image $g'$ of a geodesic curve $g$ into the surface $\mathcal{S}_{A_0A_1A_2}$
by fibre projection. $g'$ is given  by the following requirements:
\begin{enumerate}
\item[a.] First we assume that the point $P_1$ is inner point of the geodesic segment $A_iA_j$ and $P_2$ is inner point of the geodesic segment 
$A_iA_k$ $(i,j,k \in\{0,1,2\},$ $i \ne j,k, j \ne k)$.
Furthermore, we assume that there can be no midpoints at the same time. 
In this case their connecting {\it line} $\mathcal{G}^{\mathcal{S}_{A_0A_1A_2}}_{P_1P_2}$ is the fibre projected image $g'$ of the geodesic curve 
$g$ into the surface $\mathcal{S}_{A_0A_1A_2}$ where
$g$ is determined by the points $P_1^*=\mathcal{P}(P_1)$, $P_2^*=\mathcal{P}(P_2)$ and $P_3^* \in \mathcal{P}(g_{A_jA_k})$. $P_3^*$ is given by the generalized simple ratio 
(see Lemma 2.11 and Corollary 2.12). This is hereinafter referred to as the Menelaus' condition:
\begin{equation}
s^N(A_j,P_3,A_k) =  -\frac{d(A_j,P_3)}{d(P_3,A_k)}=\frac{-1}{s^N(A_j,P_1,A_i)\cdot s^N(A_i,P_2,A_k)}. \tag{4.2}
\end{equation}  
($P_3$ does not lie between the two points $A_j$ and $A_k$). 
The points $P_1^*,P_2^*,P_3^*$ determine a circle or a straight  
line in the $[x,y]$ base plane and so the corresponding surface $\mathcal{F}_g$ (cylinder or a plane) of which surface contains the geodesic curve $g$, $g'$ is 
$\mathcal{S}_{A_0A_1A_2} \cap \mathcal{F}_g$. 
\begin{rmrk}
In the previous point in formula (4.2), we chose the traditional constant $-1$ for the Menelaus' condition, 
but here another negative real number may be a suitable choice, within certain limits. 
\end{rmrk}
\item[b.] If $P_1$ and $P_2$ are the midpoints of the geodesic segments $g_{A_iA_j}$ and $g_{A_jA_k}$, respectively then their connecting {\it line} $\mathcal{G}^{\mathcal{S}_{A_0A_1A_2}}_{P_1P_2}$ 
is the fibre projected image $g'$ of the geodesic curve 
$g$ into the surface $\mathcal{S}_{A_0A_1A_2}$ where the $\theta_g$ parameter of the geodesic line $g$ is equal to the $\theta_{g_{AjA_k}}$ parameter of the geodesic line $g_{A_jA_k}$.
\item[c.] In that case if $P_1$ or $P_2$ or both points are inner points of the surface $\mathcal{S}_{A_0A_1A_2}$ then their connecting {\it line} is derived similarly to $2.a.$ configuration.
\item[d.] If $P_1$ or $P_2$ coincides with a vertex of the geodesic triangle $A_0A_1A_2$, e.g. $P_1=A_i$ and $P_2$ is an inner point of the geodesic triangle $A_0A_1A_2$ or 
lies on the geodesic line $g_{A_jA_k}$ then the
{\it line} (connecting curve) $\mathcal{G}^{\mathcal{S}_{A_0A_1A_2}}_{P_1P_2}$ is the image $g'$ of a geodesic curve $g$ into the surface $\mathcal{S}_{A_0A_1A_2}$
by fibre projection. $g'$ is derived by the following steps:
\begin{enumerate}
{\rm
\item It is a natural requirement that the requirements described in 2.a.b.c be true for new sub-triangles formed by former geodesic sides and by new surface {\it lines} 
in the triangle (see Fig.~3). Let $T$ be an inner point not on any side of the surface geodesic triangle $A_0A_1A_2$ ($T \in \mathcal{S}_{A_0A_1A_2}$). 
The surface lines $\mathcal{G}^{\mathcal{S}_{A_0A_1A_2}}_{A_0T}$, $\mathcal{G}^{\mathcal{S}_{A_0A_1A_2}}_{A_1T}$, $\mathcal{G}^{\mathcal{S}_{A_0A_1A_2}}_{A_2T}$ 
have to satisfy, that there are uniquely common inner points with the opposite geodesic side segments. 
The curve $\mathcal{G}^{\mathcal{S}_{A_0A_1A_2}}_{A_0T}$ and $g_{A_1A_2}$ meet in $P_{12}$, $\mathcal{G}^{\mathcal{S}_{A_0A_1A_2}}_{A_1T}$ and $g_{A_0A_2}$ in $P_{02}$ 
and $\mathcal{G}^{\mathcal{S}_{A_0A_1A_2}}_{A_2T}$ 
and $g_{A_0A_1}$ in $P_{12}$ lying on the geodesic segment 
$g_{A_1A_2}$ and a point $P_{02}$ on the geodesic segment $g_{A_0A_2}$. Therefore, the location of the inner point $T$ is 
uniquely determined by surface lines e.g. $T=\mathcal{G}^{\mathcal{S}_{A_0A_1A_2}}_{A_0P_{12}} \cap \mathcal{G}^{\mathcal{S}_{A_0A_1A_2}}_{A_1P_{02}}$ 
that can be given by the simple ratios $\delta_1:=s^N(A_1,P_{12},A_2)=\frac{d(A_1,P_{12})}{d(P_{12},A_2)}$ and
$\delta_2=s^N(A_2,P_{02},A_0)=\frac{d(A_2,P_{02})}{d(P_{02},A_0)}$ where $d$ is the distance function of $\NIL$ geometry.     

\item The projected images $\overset{\LARGE\frown}{A_0^*P_{12}^*}$, $\overset{\LARGE\frown}{A_1^*P_{02}^*}$ and $\overset{\LARGE\frown}{A_2^*P_{01}^*}$ 
of the above curves $\mathcal{G}^{\mathcal{S}_{A_0A_1A_2}}_{A_0P_{12}}$, $\mathcal{G}^{\mathcal{S}_{A_0A_1A_2}}_{A_1P_{02}}$ and $\mathcal{G}^{\mathcal{S}_{A_0A_1A_2}}_{A_2P_{01}}$
are circle arcs 
or line segments in the ``base plane" of $\NIL$, similarly to the projected images of the 
geodesic segments (sides of geodesic triangle $A_0A_1A_2$) $g_{A_0A_1}$, $g_{A_1A_2}$, $g_{A_2A_0}$. Moreover, $T^*$ is the common point of circle arcs or/and line segments 
$\overset{\LARGE\frown}{A_0^*P_{12}^*}$, $\overset{\LARGE\frown}{A_1^*P_{02}^*}$ and $\overset{\LARGE\frown}{A_2^*P_{01}^*}$. We note here, the simple ratio on a fibre projected surface 
{\it line} (connecting curve) $g'=\mathcal{G}^{\mathcal{S}_{A_0A_1A_2}}_{P_1P_2}$ of geodesic curve $g$ is defined by the corresponding 
simple ratio determined on the geodesic curve $g$.

Applying the Menelaus' condition to the sub-triangles we obtain, that the Ceva's theorem must be satisfied (see Theorem 4.11)
i.e. $$s^N(A_1,P_{12},A_2) \cdot s^N(A_2,P_{02},A_0)\cdot s^N(A_0,P_{01},A_1)=1.$$ Using the Lemma 4.8 follows that the corresponding Ceva's theorem true for the 
projected configuration, too i.e. for the triangle $A_0^*A_1^*A_2^*$ and the points $T^*$, $P_{01}^*$, $P_{12}^*$, $P_{02}^*$ (see Theorem 4.12). 
Moreover, other consequence of the Menelaus' condition
and the Lemma 4.8 that the simple ratios of the (projected) circle arcs or/and line segments 
\begin{equation}
\begin{gathered}
s^c(A_0^*,T^*,P_{12}^*):=\frac{\overset{\LARGE\frown}{A_0^*T^*}}{\overset{\LARGE\frown}{T^*P_{12}^*}},
s^c(A_1^*,T^*,P_{02}^*):=\frac{\overset{\LARGE\frown}{A_1^*T^*}}{\overset{\LARGE\frown}{T^*P_{02}^*}}, \\
s^c(A_2^*,T^*,P_{01}^*):=\frac{\overset{\LARGE\frown}{A_2^*T^*}}{\overset{\LARGE\frown}{T^*P_{01}^*}}, 
\end{gathered} \notag
\end{equation}
are also determined by the ratios $\delta_1$ and $\delta_2$. 
\item Therefore, location of the points $T^*$, $P_{01}^*$, $P_{12}^*$, $P_{02}^*$ uniquely determined thus the  
circle arcs or line segments $\overset{\LARGE\frown}{A_0^*T^*P_{12}^*}$, $\overset{\LARGE\frown}{A_1^*T^*P_{02}^*}$, $\overset{\LARGE\frown}{A_2^*T^*P_{01}^*}$ and by the 
Lemma 4.8. the surface lines $\mathcal{G}^{\mathcal{S}_{A_0A_1A_2}}_{A_0P_{12}}$, $\mathcal{G}^{\mathcal{S}_{A_0A_1A_2}}_{A_1P_{02}}$, $\mathcal{G}^{\mathcal{S}_{A_0A_1A_2}}_{A_2P_{01}}$
are given, respectively.}
\end{enumerate}
\end{enumerate}
\begin{figure}[ht]
\centering
\includegraphics[width=12cm]{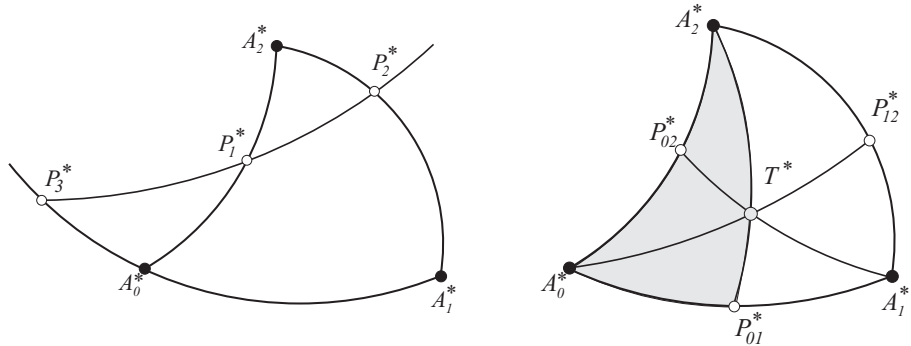}
\caption{Projected image of surface {\it line} (connecting curve) $\mathcal{G}^{\mathcal{S}_{A_0A_1A_2}}_{P_1P_2}$ and the projected image of a Ceva's configuration. 
}
\label{}
\end{figure}
\end{definition}
\begin{corollary}
Based on all this, it can be seen that Menelaus' theorem does not follow from the structure of $\NIL$ geometry. 
However, as can be seen above, the Menelaus' condition plays an important role in defining {\it lines} on surfaces by geodesic lines given triangle.
\end{corollary}
Using the above Menelaus' condition, similar to the Euclidean proof, we obtain the $\NIL$ Cava's theorem:
\begin{Theorem}
If $T$ is a point not on any side of a geodesic triangle $A_0A_1A_2$ in $\NIL$ space such that
the curves $\mathcal{G}^{\mathcal{S}_{A_0A_1A_2}}_{A_0T}$ and $g_{A_1A_2}$ meet in $P_{12}$, $\mathcal{G}^{\mathcal{S}_{A_0A_1A_2}}_{A_1T}$ and $g_{A_0A_2}$ in $P_{02}$, 
and $\mathcal{G}^{\mathcal{S}_{A_0A_1A_2}}_{A_2T}$ and $g_{A_0A_1}$ in $P_{01}$, 
then $$s^N(A_0,P_{01},A_1)s^N(A_1,P_{12},A_2)s^N(A_2,P_{02},A_0) = 1.$$ ~ ~ ~ $\square$
\end{Theorem}
Using the Lemma 4.8 follows that the corresponding Ceva's theorem true for the 
{\it projected configuration} too i.e. for the triangle $A_0^*A_1^*A_2^*$ and the points $T^*$, $P_{01}^*$, $P_{12}^*$, $P_{02}^*$. 
\begin{Theorem}
If $T^*$ is a point not on any side of {\it circle arc} triangle (the projected image of a geodesic triangle in general type) $A_0^*A_1^*A_2^*$ in the base plane of the $\NIL$ space such that
the arcs (or line segments)  $\overset{\LARGE\frown}{A_0^*T^*}$ and $\overset{\LARGE\frown}{A_1^*A_2^*}$ meet in $P_{12}^*$, $\overset{\LARGE\frown}{A_1^*T^*}$ and 
$\overset{\LARGE\frown}{A_0^*A_2^*}$ in $P_{02}^*$, and $\overset{\LARGE\frown}{A_2^*T^*}$ and $\overset{\LARGE\frown}{A_0^*A_1^*}$ in $P_{01}^*$, 
then $$s^c(A_0^*,P_{01}^*,A_1^*)s^c(A_1^*,P_{12}^*,A_2^*)s^c(A_2^*,P_{02}^*,A_0^*) = 1.$$ ~ ~ ~ $\square$
\end{Theorem}
\begin{rmrk}
Using the previous notions and theorems, similar to Euclidean concepts, we can define, for example, the circumscribed 
circle of a geodesic triangle and their centres, the centroid of a geodesic triangle as the point where the three medians of the triangle meet. 
(a median of a geodesic triangle $A_0A_1A_2$ in the $\NIL$ space is a surface {\it line} $\subset \mathcal{S}_{A_0A_1A_2})$ from one vertex to the mid point on the 
opposite side of the triangle). But we will examine these in a forthcoming paper. 
\end{rmrk}

Similar problems in other homogeneous Thurston geometries
represent another huge class of open mathematical problems. For
$\SOL$, $\SLR$ geometries only very few results are known
\cite{CsSz16}, \cite{MSz}, \cite{MSz12}, \cite{Sz13-1}, \cite{Sz13-2}, \cite{Sz12-1}, \cite{Sz19}, \cite{Sz11-1}, \cite{VSz19}. Detailed studies are the objective of
ongoing research.
\medbreak

\end{document}